\documentclass[preprint,12pt]{elsarticle}
\usepackage{amssymb}
\usepackage{xcolor}
\usepackage{float}
\usepackage{amsmath}
\usepackage{amsthm}
\newtheorem{theorem}{Theorem}[section]
\begin{document}

\begin{frontmatter}

\title{Convergence of a Ramshaw-Mesina Iteration} %% 

\author[label1]{Aytekin \c{C}{\i}b{\i}k} %% Author name
\author[label2]{ and William Layton}
\affiliation[label1]{organization={Department of Mathematics, Gazi University},%Department and Organization 
            city={Ankara},
            postcode={06560}, 
            country={Turkey}}
\affiliation[label2]{organization={Department of Mathematics, University of Pitsburgh},%Department and  
            city={Pittsburgh},
            postcode={15260}, 
            state={PA},
            country={USA}}
%% Abstract
\begin{abstract}
In 1991 Ramshaw and Mesina introduced a clever synthesis of penalty methods
and artificial compression methods. Its form makes it an interesting option to
replace the pressure update in the Uzawa iteration. The result, for the Stokes
problem, is%
\begin{equation}\label{system}
\left\{
\begin{array}
[c]{cc}%
Step\text{ }1: & -\triangle u^{n+1}+\nabla p^{n}=f(x),\text{ in }\Omega,\text{
}u^{n+1}|_{\partial\Omega}=0,\\
Step\text{ }2: & p^{n+1}-p^{n}+\beta\nabla\cdot(u^{n+1}-u^{n})+\alpha
^{2}\nabla\cdot u^{n+1}=0.
\end{array}
\right.
\end{equation}
For saddle point problems, including Stokes, this iteration converges under a condition similar to the one required for Uzawa iteration.

\end{abstract}



%% Keywords

\end{frontmatter}



\section{Introduction}

In 1991 Ramshaw and Mesina \cite{RM91} gave a clever synthesis of
regularization of incompressibility ($\nabla\cdot u=0$) by artificial
compression ($p_{t}+\alpha^{2}\nabla\cdot u=0$) and a reformulation of
penalization ($\frac{\partial}{\partial t}\left(  p+\beta\nabla\cdot
u=0\right)  $). For $\alpha,\beta$ suitably balanced, their regularization%
\begin{equation}
p_{t}+\beta\nabla\cdot u_{t}+\alpha^{2}\nabla\cdot u=0\label{eq:RM}%
\end{equation}
sped up satisfaction of incompressibility by several orders of magnitude.
Since (\ref{eq:RM}) reduces to $\nabla\cdot u=0$\ at steady state
($p_{t}=0,u_{t}=0$), it is natural to consider replacing the pressure update
step in the Uzawa iteration by a discrete form of (\ref{eq:RM}). For the
Stokes problem, this yields: given $u^{n},p^{n}$%
\begin{equation}
\left\{
\begin{array}
[c]{ccc}%
Step\text{ }1: &  & -\triangle u^{n+1}+\nabla p^{n}=f(x),\text{ in }%
\Omega,\text{ }u^{n+1}|_{\partial\Omega}=0,\\
&  & \\
Step\text{ }2: &  & p^{n+1}-p^{n}+\beta\nabla\cdot(u^{n+1}-u^{n})+\alpha
^{2}\nabla\cdot u^{n+1}=0.
\end{array}
\right.  \label{eqRMIteration}%
\end{equation}
If $\beta=0$ this reduces to the standard Uzawa iteration. Section 2 proves
convergence for (\ref{eqRMIteration}) under the condition $\beta+\alpha
^{2}<M^{-1}$, where $M$ is the maximum of the Rayleigh quotient for the Schur
complement. Section 3 provides a numerical test.

\textbf{Some related work}. The utility of (\ref{eq:RM}) to timestep to steady
state developed in directions complemented herein by Ramshaw and Mousseau
\cite{RM91b} and McHugh and Ramshaw \cite{MR95}. The (short) proof herein
builds on the analysis of the Uzawa iteration in Bacuta \cite{B06}. The
regularization (\ref{eq:RM}) was explored in different directions in
\cite{RM91b},\cite{RM94} \cite{duk93}, \cite{MR95}. At this point, experience
with (\ref{eqRMIteration}) is insufficient to classify its advantages,
disadvantages, complexity differences. There are many possible further
developments paralleling those of the Uzawa iteration in e.g. \cite{BWY89},
\cite{EG94}, (among hundreds of papers).
\section{Proof of convergence}
We follow the framework of Bacuta \cite{B06} which we now summarize.  Let
$\mathbf{V},P$ denote two Hilbert spaces.  Denote linear operators $A:\mathbf{V}\rightarrow\mathbf{V}^{\ast}$ and
$B:V\rightarrow P,$ so that $B^{\ast}:P\rightarrow\mathbf{V}^{\ast}$. In this
framework, we consider the reformulation of (\ref{eqRMIteration}) as%
\begin{equation}
\left\{
\begin{array}
[c]{c}%
Au^{n+1}+B^{\ast}p^{n}=f,\\
p^{n+1}-p^{n}-\beta B(u^{n+1}-u^{n})-\alpha^{2}Bu^{n+1}=0.
\end{array}
\right. \label{eq:Reformulation}%
\end{equation}
For the Stokes problem this corresponds to $A=-\triangle,B^{\ast}%
=\nabla,B=-\nabla\cdot,\mathbf{V}=(H_{0}^{1}(\Omega))^d,\, P= L_{0}^{2}(\Omega)$, with $f\in\mathbf{V}^{\ast}$. The Schur complement operator $\mathcal{A}$ is
assumed invertible and that%
\[
\mathcal{A}:=BA^{-1}B^{\ast}:P\rightarrow P
\]
 is a self-adjoint, bounded, positive definite operator. These properties are proven under mild conditions in Bacuta \cite{B06}, Lemma 2.1 p.2635. In particular, Bacuta
\cite{B06} shows that there are positve constants $m,M$ such that, for
$(\cdot,\cdot),|\cdot|$ the \ $P$-inner product and norm,%
\begin{equation}
0<m|q|^{2}\leq(\mathcal{A}q,q)\leq M|q|^{2}\text{ for all }0\neq q\in
P.\label{eq:KeyEstimate}%
\end{equation}
Under these conditions, the $\mathcal{A}$-inner product and norm,
$(p,q)_{\mathcal{A}}:=(\mathcal{A}p,q),|p|_{\mathcal{A}}:=(p,p)_{\mathcal{A}%
}^{1/2},$ are well-defined. We prove the following.

\begin{theorem}
Under the above assumptions, in particular (\ref{eq:KeyEstimate}), the
iteration (\ref{eqRMIteration}) converges if $\beta\geq0$ and\ $\beta
+\alpha^{2}<\frac{1}{M}.$
\end{theorem}

\begin{proof} Step 1 of (\ref{eq:Reformulation}) is used to eliminate the
velocity in a standard way using $u^{n+1}=A^{-1}(f-B^{\ast}p^{n})$ which
reduces Step\ 2 to%
\[
p^{n+1}-p^{n}+\beta\mathcal{A}(p^{n}-p^{n-1})+\alpha^{2}\mathcal{A}%
p^{n}=\alpha^{2}BA^{-1}f.
\]
The true $p$ satisfies this exactly. Thus the error $e^{n}=p-p^{n}$ satisfies,
by subtraction and rearrangement,%
\[
e^{n+1}-\left[  I-(\beta+\alpha^{2})\mathcal{A}\right]  e^{n}-\beta
\mathcal{A}e^{n-1}=0.
\]
Take the inner-product with $e^{n+1}.$ This gives:
\[
|e^{n+1}|^{2}-\mathcal{(}e^{n},e^{n+1})+(\beta+\alpha^{2})\mathcal{(}%
e^{n},e^{n+1})_{\mathcal{A}}-\beta\mathcal{(}e^{n-1},e^{n+1})_{\mathcal{A}}=0.
\]
\ Using the polarization identity on $\mathcal{(}e^{n},e^{n+1}),\mathcal{(}%
e^{n},e^{n+1})_{\mathcal{A}}$\ and $\mathcal{(}e^{n-1},e^{n+1})_{\mathcal{A}}$
gives%
\begin{gather*}
|e^{n+1}|^{2}-\left[  \frac{1}{2}|e^{n}|^{2}+\frac{1}{2}|e^{n+1}|^{2}-\frac
{1}{2}|e^{n+1}-e^{n}|^{2}\right]  +\\
(\beta+\alpha^{2})\left[  \frac{1}{2}|e^{n}|_{\mathcal{A}}^{2}+\frac{1}%
{2}|e^{n+1}|_{\mathcal{A}}^{2}-\frac{1}{2}|e^{n+1}-e^{n}|_{\mathcal{A}}%
^{2}\right]  \\
-\beta\left[  \frac{1}{2}|e^{n-1}|_{\mathcal{A}}^{2}+\frac{1}{2}%
|e^{n+1}|_{\mathcal{A}}^{2}-\frac{1}{2}|e^{n+1}-e^{n-1}|_{\mathcal{A}}%
^{2}\right]  =0.
\end{gather*}
Multiply by 2 and regroup terms by%
\begin{gather}
\left[  |e^{n+1}|^{2}+\beta|e^{n}|_{\mathcal{A}}^{2}+\alpha^{2}|e^{n+1}%
|_{\mathcal{A}}^{2}+\left\{  |e^{n+1}-e^{n}|^{2}-(\beta+\alpha^{2}%
)|e^{n+1}-e^{n}|_{\mathcal{A}}^{2}\right\}  \right]  \label{eq:CriticalStep}\\
-\left[  |e^{n}|^{2}+\beta|e^{n-1}|_{\mathcal{A}}^{2}+\alpha^{2}%
|e^{n}|_{\mathcal{A}}^{2}+\left\{  |e^{n}-e^{n-1}|^{2}-(\beta+\alpha
^{2})|e^{n}-e^{n-1}|_{\mathcal{A}}^{2}\right\}  \right]  \nonumber\\
+\left\{  |e^{n}-e^{n-1}|^{2}-(\beta+\alpha^{2})|e^{n}-e^{n-1}|_{\mathcal{A}%
}^{2}\right\}  \nonumber\\
+2\alpha^{2}|e^{n}|_{\mathcal{A}}^{2}+\beta|e^{n+1}-e^{n-1}|_{\mathcal{A}}%
^{2}=0.\nonumber
\end{gather}
Notice that in three places (in braces, \{$\cdot$\}) there occurs
$|q|^{2}-(\beta+\alpha^{2})|q|_{\mathcal{A}}^{2}$ with $q=e^{n+1}-e^{n}$ and
$q=e^{n}-e^{n-1}$. Using $0<m|q|^{2}\leq(\mathcal{A}q,q)\leq M|q|^{2}$ we
bound this term as%
\[
\left[  1-M(\beta+\alpha^{2})\right]  |q|^{2}\leq|q|^{2}-(\beta+\alpha
^{2})|q|_{\mathcal{A}}^{2}\leq\left[  1-m(\beta+\alpha^{2})\right]  |q|^{2}.
\]
Thus, this term is positive if $\beta+\alpha^{2}<1/M.$ Thus,
(\ref{eq:CriticalStep}) takes the abstract form $\left[  E^{n+1}\right]
-\left[  E^{n}\right]  +P^{n}=0$ where both $E\&P$ are non-negative. Summing
we have $E^{N}$ and $\sum_{n=1}^{N}P^{n}$, a series with non-negative terms,
are uniformly bounded. Hence $\sum_{n=1}^{\infty}P^{n}$\ is convergent and the
$n^{th}$ term $P^{n}\rightarrow0$ as $n\rightarrow\infty$ where%
\[
P^{n}=\left\{  |e^{n}-e^{n-1}|^{2}-(\beta+\alpha^{2})|e^{n}-e^{n-1}%
|_{\mathcal{A}}^{2}\right\}  +2\alpha^{2}|e^{n}|_{\mathcal{A}}^{2}%
+\beta|e^{n+1}-e^{n-1}|_{\mathcal{A}}^{2}.
\]
This implies $e^{n}\rightarrow0$, concluding the proof.

\end{proof}
\section{Numerical Test}
We compare (\ref{system}) to the standard Uzawa algorithm for the lid-driven cavity. The domain is a unit square with the top lid sliding in the positive $x$-direction according to $u(x,1)=(g(x),0)^{T},\, g(x)=4x(1-x) $. The boundary conditions are $u=0$ on the other walls. FreeFem++ is used for computations and the inf-sup stable Taylor-Hood finite element pair is used for the test.\par
The stopping criteria for the iterations is:
\begin{align}\label{stop}
	\max\{\|u^{n+1}-u^{n}\|,\|p^{n+1}-p^{n}\|\}\leq 10^{-6}.
\end{align}
According to \cite{bans2003}, the value of $\alpha^2$ should be $\alpha^2<2$. Our numerical tests confirmed that for any $\alpha^2\geq2$ the method diverges and the optimal selection of $\alpha^2$ is $\alpha^2=1.5$ independent of $\beta$.
\begin{table}[H]
	\centering
	\begin{tabular}{|c|c|c|c|c|c|c|c||}
		\hline
		Mesh  & $\beta=0$ &$\beta=10^{-4}$& $\beta=10^{-2}$ &$\beta=0.1$& $\beta=0.2$  \\
		\hline
		$10\times10$&74&74&74&75&89\\
		\hline
		$20\times20$&75&75&75&77&85
		\\
		\hline
		$40\times40$&77&77&77&78&85
		\\
		\hline
	\end{tabular}
	\caption{ Iterations  for  different meshes for varying $\beta$ and $\alpha^2=1.5$.}
	\label{tab:iter}
\end{table}
Table \ref{tab:iter} indicates that mesh refinement has a very slight effect on convergence and increasing $\beta$ up to $0.1$ has almost no effect on number of iterations. For $\beta> 0.2$, (\ref{system}) diverged due to the violation of convergence criteria. Further tests for smaller $\alpha$ and larger $\beta$ did not show a significant change in number of iterations. In this test, the method had almost the same convergence properties as standard Uzawa.\par
Next we also tested the effect of solution regularity by changing the lid to $u(x,1)=(1,0)^{T}$ with fixed $\alpha$ and $\beta$ values ($\alpha^2=1.5$, $\beta= 0.05$)
\begin{figure}[H]
	\centering
	\includegraphics[width=0.7\linewidth]{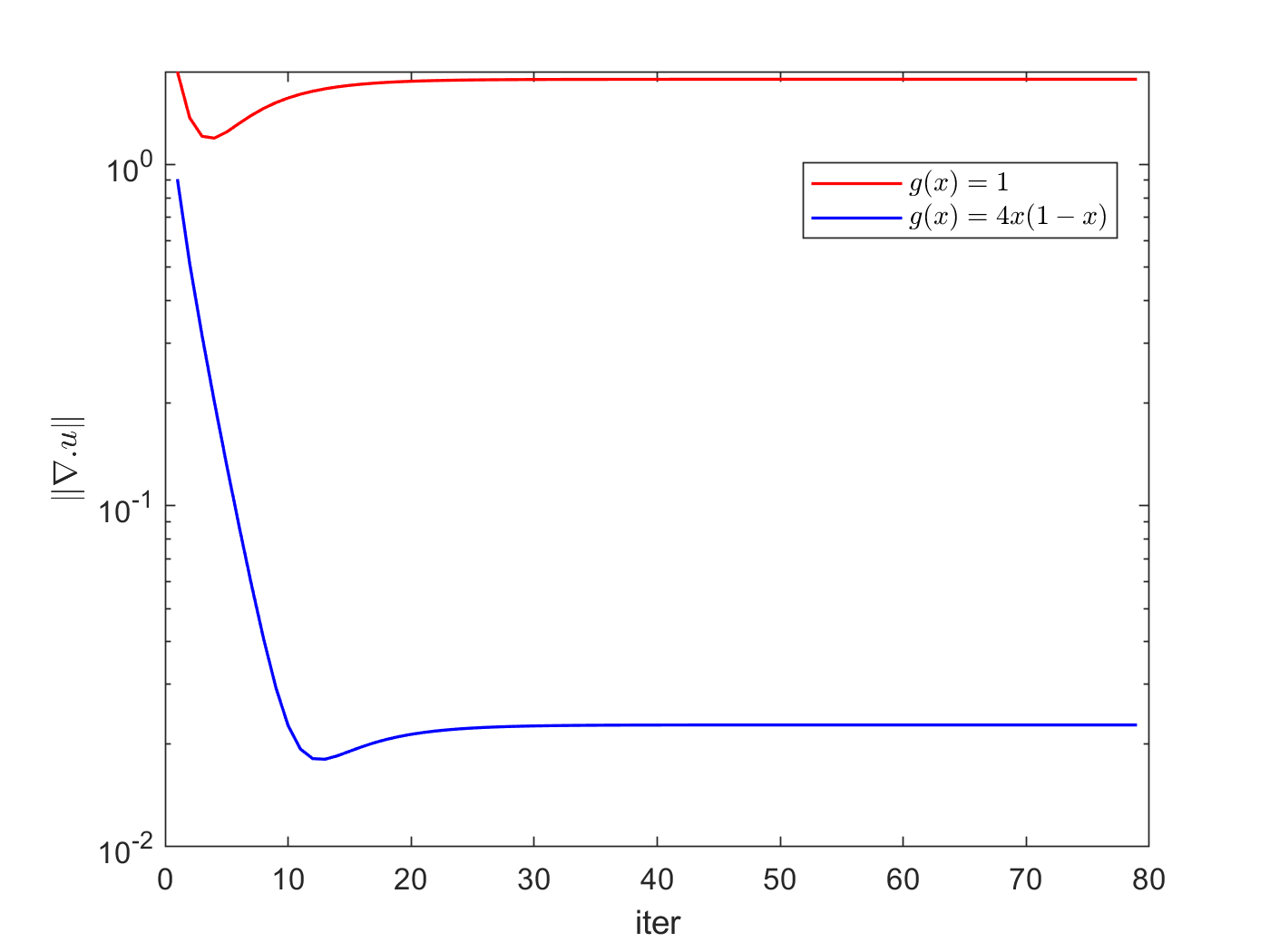}
	\caption{Comparison of $\|\nabla . u\|$ for different selections of $g(x)$.}
	\label{fig:divergence}
\end{figure}
Figure \ref{fig:divergence} indicates that solution regularity has a notable effect on $\|\nabla . u\|$.
\section{Conclusion}
If the pressure of Uzawa is replaced as in (\ref{system}), convergence seems to be similar to standard Uzawa. This suggests the dramatic improvement observed by Ramshaw and Mesina is linked to explicit time stepping. Replacing Step 1 in (\ref{system}) by a first order Richardson step is therefore where the impact of Step 2 should be next explored.\\

\noindent {\bf Funding} The research was partially supported by NSF grant DMS-2110379 and TUBITAK grant BIDEB2219.

\end{document}